\input amstex
\documentstyle{amsppt}
\input bull-ppt
\keyedby{bull399e/PAZ}

\define\CPtwo{\Bbb C\Bbb P^2}
\define\CPtwobar{\overline{\CPtwo}}
\define\SU{\mathop{SU}\nolimits}
\define\SO{\mathop{SO}\nolimits}


\topmatter
\cvol{29}
\cvolyear{1993}
\cmonth{July}
\cyear{1993}
\cvolno{1}
\cpgs{63-69}
\title
The genus-minimizing property\\
of algebraic curves
\endtitle
\shorttitle{Genus-minimizing Curves}
\author	P. B. Kronheimer
\endauthor
\address Mathematical Institute, 
24--29 St.~Giles, Oxford OX1 3LB
\endaddress
\ml kronheim\@maths.oxford.ac.uk
\endml
\date November 17, 1992\enddate
\abstract A viable and still unproved conjecture states 
that, if $X$
is a smooth algebraic surface and $C$ is a smooth algebraic
curve in $X$, then $C$ realizes the smallest possible genus
amongst all smoothly embedded $2$-manifolds in its homology
class. A proof is announced here for this conjecture, for a
large class of surfaces $X$, under the assumption that the
normal bundle of $C$ has positive degree.
\endabstract
\endtopmatter


\document

\heading 1. Introduction
\endheading
If $X$ is a smooth $4$-manifold and $\xi$ is a 
$2$-dimensional homology
class in $X$, one can always represent $\xi$ geometrically 
by an oriented
$2$-dimensional surface $\varSigma$, smoothly embedded in 
the $4$-manifold.
Depending on $X$ and $\xi$ however, the genus of 
$\varSigma$ may have to be
quite large: it is not always possible to represent $\xi$ 
by an embedded
sphere.  It is natural to ask for a representative whose 
genus is as small
as possible, or at least to enquire what the genus of such 
a minimal
representative would be.  Although not much is known about 
this question in
general, there is an attractive conjecture concerning the 
case that $X$ is
the manifold underlying a smooth complex-algebraic 
surface.  The conjecture
is best known in the case that $X$ is the complex 
projective plane $\CPtwo$,
in which case it is often attributed to Thom, but the 
statement seems still 
to be viable more generally \cite{1}.

\thm{Conjecture 1}
Let $X$ be a smooth algebraic surface and $\xi$ a homology 
class
carried by a smooth algebraic curve $C$ in $X$.  Then $C$ 
realizes the
smallest possible genus amongst all smoothly embedded  
$2$-manifolds
representing $\xi$.
\ethm

The attractiveness of this conjecture stems from the 
connection
to which it points, 
between low-dimensional topology and complex geometry.  
Through the work
of Donaldson in particular, this connection is now a 
familiar feature of
differential topology in dimension~$4$, and the techniques 
of gauge theory
provide a natural starting point for an approach to the 
problem.  The
conjecture was proved in \cite{5} for the special case 
that $X$ is a $K3$
surface, and the results of that paper also gave a lower 
bound for the genus
of an embedded surface in more general complex surfaces. 
However, in all
applicable cases other than $K3$, the lower bound proved 
in \cite{5} falls
short of Conjecture~1. The purpose of this present paper 
is to describe a
result which establishes the correctness of the conjecture 
for a large class
of complex surfaces.  The hypotheses of Theorem~2 still 
exclude the
tantalizing case of $\CPtwo$, but conditions (a) and (c) 
of the theorem
admit very many (and conjecturally all) simply connected 
surfaces $X$ of
general type and odd geometric genus.  Full details of the 
proof will appear
later in \cite{4}.

\proclaim {Theorem 2}
The above conjecture holds at least under the following 
assumptions
concerning $X$ and $C$\RM:
\par
{\rm(a)} the surface $X$ is simply connected\RM;
\par
{\rm(b)} the self-intersection number $C\cdot C$ is 
positive\,\RM;
\par
{\rm(c)} there is a class $\omega \in H_2(X,\Bbb C)$ dual 
to a
holomorphic $2$-form on $X$\RM, such that $q_k(\omega +
\overline \omega) > 0$ for sufficiently large $k$\RM, 
where $q_k$
		denotes Donaldson's polynomial invariant.
\endproclaim

Some comments are needed concerning the third hypothesis.  
Donaldson's {\it
polynomial invariants\/} \cite{2} are homogeneous 
polynomial functions $q_k$
on $H_2(X)$, defined using instanton moduli spaces for 
structure group
$\SU(2)$; their degree depends on the parameter $k$ as 
well as the homotopy
type of $X$.  Condition (c) appears in \cite{8}, where it 
is shown that this
condition will hold for a surface $X$ of general type 
provided that:
\roster
  \item"(i)"	the geometric genus $p_g(X)$ is odd; and

  \item"(ii)"	the canonical linear system of the minimal 
model of $X$ 
		contains a smooth curve.
\endroster
(This result also rests on some more technical material in 
\cite{6}.)  The
first of these two conditions ensures that the polynomial 
invariants have even
degree and certainly is essential as long as one considers 
the invariants
associated with the structure group $\SU(2)$, though the 
$\SO(3)$ moduli
spaces can be used to treat some of the remaining cases.  
The importance of
condition (ii) is less clear, but it does indicate that 
one should expect
(c) to be a rather general property of complex surfaces 
whenever the
polynomial invariants have even degree. 

The basic material of the proof of Theorem~2 is the same 
as that of the main
theorem of \cite{5}, namely, the moduli spaces of 
instantons on $X$ having a
singularity along an embedded surface.  The structure of 
the argument,
however, is
rather different.  The difficulty of embedding 
$2$-dimensional
surfaces in four dimensions stems from one's inability to 
remove unwanted
self-intersection points of an {\it immersed\/} surface, 
even when these
intersections cancel algebraically in plus--minus pairs; 
this is the failure
of the Whitney lemma in dimension $4$ and lies at the 
heart of $4$-manifold
topology and all its problems.  The first stage of the 
proof of Theorem~2 is
the construction of invariants which measure an 
obstruction to the removal
of such pairs of intersection points. Given an immersed 
surface $\varSigma$
with normal crossings in a $4$-manifold $X$, we use the 
moduli spaces of
singular instantons to define an invariant of the pair 
$(X,\varSigma)$; this
will be a distinguished function
$ 
R_d(s) : H_2(X) \to \Bbb R$, 
taking the form of a homogeneous polynomial of degree $d$ 
on $H_2(X)$ and
a finite Laurent series in the formal variable $s$. This 
invariant has
the property that the order of vanishing of $R_d$ at the 
point $s=1$ gives
an upper bound on the number of positive-signed 
intersection points in
$\varSigma$ which can ever be removed by a homotopy of the 
immersion.  The
second stage of the proof is to show that, in the case of 
an algebraic curve
in a suitable algebraic surface, the invariants $R_d(s)$ 
give information
which is sharp enough to establish the assertion of 
Conjecture~1; this
entails proving a nonvanishing theorem for the value of 
$R_d(s)$ at $s=1$.
Using ideas from \cite{8}, we shall in fact show that, if 
$C'$ is an
irreducible algebraic curve with a single ordinary double 
point in a
suitable complex surface $X$, then $R_d(1)$ is positive 
for the pair
$(X,C')$ when evaluated on a class $\omega + 
\overline\omega$ as in 
(c), so
showing that the self-intersection point in $C'$ cannot be 
removed by any
homotopy. Theorem~2 is easily deduced from this.

The structure of this proof is closely modeled on 
Donaldson's proof of the
indecomposability of complex surfaces in \cite{2}, but a 
closer parallel
still is in \cite{1}, where a similar strategy was used in 
connection with
Conjecture~1.  In that paper, use was made of the 
instanton moduli spaces
associated with a branched double cover $\widetilde X \to 
X$, branched along
$\varSigma$. The moduli spaces of singular instantons 
which we use here can, in
some circumstances, be interpreted as moduli spaces of 
suitably {\it
equivariant\/} instanton connections on such a covering 
manifold
(equivariant, that is, under the covering involution). The 
replacing of the
full moduli space on $\widetilde X$ by this equivariant 
part can be seen as the
main difference in the framework of the argument between 
\cite{1} and the
present paper.  The other main new ingredient here is the 
organization of
the invariants obtained from the singular instantons to 
form the finite
Laurent series $R_d$.  The remaining two sections of this 
paper provide some
further details of the proof; some of this material is 
rather technical,
particularly in \S3. It seems likely that a change of 
strategy in this
last part of the argument will eventually lead to a 
slightly more general
result than Theorem~2. 


\heading 2. Obstructions to removing intersection points
\endheading
Let $X$ be a smooth, oriented, simply connected closed 
$4$-manifold and
$\varSigma$ an embedded (rather than just immersed) 
orientable surface in
$X$.  Given a Riemannian metric on $X$, it was shown in 
\cite{5} how one
can construct moduli spaces $M^\alpha_{k,l}$ associated to 
the pair; roughly
speaking, $M^\alpha_{k,l}$ parametrizes the finite-action 
anti-self-dual
$\SU(2)$ connections on $X\setminus\varSigma$ with the 
property that, near
to $\varSigma$, the holonomy around small loops linking 
the surface is
asymptotically 
$\exp 2\pi i(\stack{-\alpha}{0}\ \stack{0}{\alpha})$.
Here $\alpha$ is a real parameter in the interval $(0, 
1/2)$ and $k$ and $l$
are the integer topological invariants of such 
connections: the ``instanton''
and ``monopole'' numbers.  For a generic choice of 
Riemannian metric and away
from the flat or reducible connections in the moduli 
space, $M^\alpha_{k,l}$
is a smooth manifold of dimension
$$
		8k + 4l - 3(b^+ + 1) - (2g-2) \tag1
$$
where $b^+$ is the dimension of a maximal positive 
subspace for the
intersection form on $H^2(X)$ and $g$ is the genus of 
$\varSigma$.   In the
case that $b^+$ is odd, the dimension is even, and we 
write it as $2d(k,l)$,
so that $d = d(k,l)$ is half of (1). Following Donaldson's 
definition
\cite{2} of the polynomial invariants $q_k$, it was shown 
in \cite{5} that
the moduli spaces $M^\alpha_{k,l}$ can be used, when $b^+$ 
is odd, to define
a homogeneous polynomial function of degree $d(k,l)$,
$$
		q_{k,l} : H_2(X) \to \Bbb R.
$$
(In \cite{5}, this polynomial was defined only on the 
orthogonal complement
of $[\varSigma]$ in $H_2(X)$, and this is all we will 
actually need; the
definition, however, can be extended to the whole homology 
group.  Also, a
`homology orientation' of $X$ is needed to fix the overall 
sign.)
When $b^+$ is at least three, the polynomial $q_{k,l}$ is 
independent of the
parameter $\alpha$ and the Riemannian metric; it is an 
invariant of the pair
$(X,\varSigma)$.

Because of the way $k$ and $l$ enter the dimension formula 
(1), the degree
of $q_{k-1, l+2}$ is the same as the degree of $q_{k,l}$. 
It is natural to
combine all the polynomials of a given degree into one 
Laurent
series: 
$$
	R_d(s)  =  \sum_{d(k,l) = d}  s^l q_{k,l}.
$$
This series is actually finite in both directions, though 
this is not
obvious a priori from its definition.  Note that, 
depending on the parity of
the genus $g$ and the value of $b^+$ mod $4$, the 
invariant $R_d$ will be
defined only for $d$ of one particular parity.  Flat 
connections on the
complement of $\varSigma$ can cause difficulties in 
defining the invariants
directly from the moduli spaces, but these can be 
overcome, for example, by
a device such as that described in \cite{7}.

Having defined invariants for embedded surfaces, we can 
now define
invariants for immersed surfaces with normal crossings.  
It would seem
feasible to do this by directly using gauge theory on the 
complement of the
immersed surface, but a short-cut is available.  We shall 
convert such an
immersed surface $\varSigma$ into an embedded surface by 
blowing up $X$ at
the intersection points.  This is the process modeled on 
the situation in
complex geometry, where a curve $\varSigma$ with a normal 
crossing at $p$ is
replaced by its proper transform, which is a smooth curve 
in the new surface
$\widetilde X = X \# 
\overline{\Bbb{C}\Bbb{P}}^2$, the blow-up of $X$ at $p$.  
In the
$C^\infty$ case, the model is just the same; there are 
really two different
cases according to the sign of the intersection point, but 
no essential
difference in the local picture.  Thus we obtain an 
embedded surface 
$\widetilde
\varSigma$ in a new manifold $\widetilde X = X \# 
n\CPtwobar$.  We
define the invariants $q_{k,l}$ and $R_d$ for 
$(X,\varSigma)$ to be the
restriction of the invariants of $(\widetilde X, 
\widetilde\varSigma)$ to $H_2(X)
\subset H_2(\widetilde X)$. 

The next stage in the argument is to see how the invariant 
$R_d(s)$ changes
when the immersion of $\varSigma$ in $X$ is changed by a 
homotopy.  During a
homotopy, double points can appear and disappear in 
$\varSigma$ in quite
complicated ways, but standard theory says that after a 
small perturbation
any such changes can be broken down into a combination of 
moves, each of
which is one of six types. One has to consider the 
following three standard
modifications and their inverses (see \cite{3} for 
pictures and explanations
of these):
\roster
  \item"(a)"	introduce a positive double point by a twist 
move;

  \item"(b)"	introduce a negative double point by a twist 
move;

  \item"(c)"	introduce a cancelling plus--minus pair by a 
finger move.
\endroster
(We should emphasize that we are talking about homotopies 
whose starting and
finishing points are immersions; only (c) can be achieved 
by a homotopy {\it
through\/} immersions.)  The change in $R_d$ under each of 
these three moves
is summarized by the following proposition.

\proclaim {Proposition 3}
Let $\varSigma$ be obtained from $\widehat\varSigma$ by 
one of the moves
just described\RM, and let $R_d$ and $\widehat R_d$ be the 
invariants for
$(X,\varSigma)$ and $(X, \widehat\varSigma)$.  Then\RM, 
according to the three
cases\RM, we have\,\RM:
\roster
  \item"(a)"	$R_d(s) = (1 - s^{-2}) \widehat R_d(s)$\RM;

  \item"(b)"	$R_d(s) = \widehat R_d(s)$\RM;

  \item"(c)"	$R_d(s) = (1 - s^{-2}) \widehat R_d(s)$.
\endroster
\endproclaim

As a consequence of these relations, the order of 
vanishing of $R_d(s)$ at
$s = 1$ increases by one every time a positive double 
point is introduced by
either of the moves (a) or (c) and decreases by one every 
time a
positive double point disappears.  So, as was stated in 
the introduction,
the order of vanishing of $R_d$ at $s=1$ puts an upper 
bound on the number
of positive double points which can be removed.

Note also that the invariants $R_d$ are unable to detect 
subtleties of
knotting: two embedded surfaces of the same genus will 
have the same
invariant if the embeddings are homotopic.  (There is also 
a simple formula
for how $R_d$ changes when the genus of $\varSigma$ is 
increased by summing
with a torus in $X$, but this involves aspects which would 
take us too far
afield; see \cite{4}.)

The proof of Proposition~3 involves some rather simple 
gluing arguments.
Consider (a) for example.  After introducing the positive 
double
point to form $\varSigma$ from $\widehat\varSigma$, the 
definition of the
invariants for immersed surfaces tells us to remove the 
double point by
blowing up, to get $(\widetilde X, \widetilde \varSigma)$.  
Examining the overall
effect, we find that, up to diffeomorphism, 
$(\widetilde X, \widetilde\varSigma)$ is
the connect sum of the pair $(X,\widehat\varSigma)$ with 
the pair
$(\CPtwobar, \overline C)$, where $C$ is a conic curve in 
the projective
plane.  One must analyse this connected sum of pairs to 
show that the
invariants for $\widetilde\varSigma$ and 
$\widehat\varSigma$ are related by 
$  q_{k,l} = \widehat q_{k,l} - \widehat q_{k-1, l+2}$.
(The sign here is rather subtle and crucial to the 
argument.)  Technical
aspects of the gluing construction, as well as some 
aspects of the algebraic
geometry in \S3, can be considerably simplified by using the
fact that the moduli spaces $M^\alpha_{k,l}$ for $\alpha = 
1/4$ are
essentially equivalent to moduli spaces of equivariant 
connections  on a
branched double cover, or to orbifold connections in the 
case that a global
double covering does not exist.  So, in the example above, 
rather than think
of forming the connect sum along a pair $(S^3, S^1)$, one 
may think instead
of gluing across a copy of $S^3$, with invariance imposed 
under an involution.


\heading 3. The complex case
\endheading
Suppose now that $X_1$ is a smooth complex surface and 
$C_1$ is a smooth
algebraic curve in $X_1$.  If the self-intersection number 
of $C_1$ is
positive, and we wish to prove Conjecture~1 for the 
homology class $[C_1]$,
then it turns out to be enough to tackle instead the 
homology class $n[C_1]$
for any large $n$ (see \cite{1} or \cite{5} for this 
elementary
construction).  So let $C_2$ be a smooth curve in the 
linear system
$|nC_1|$.  Taking $n$ large enough, we may suppose that 
the linear system
$|C_2|$ contains, in addition to this smooth curve, an 
irreducible singular
curve $C'_2$ with a single ordinary double point.  We can 
look at $C'_2$ as
an immersed $2$-manifold with a single normal crossing of 
positive sign; its
genus is one less than the genus of $C_2$.  Suppose the 
conjecture fails for
the homology class of $[C_2]$. Then we can find an {\it 
embedded\/} surface
$\varSigma$ in the this class, with the same genus as the 
immersed surface
$C_2'$.  Since $X_1$ is simply connected, $\varSigma$ and 
$C_2'$ will be
homotopic, and it follows from the results of the previous 
section that the
invariant $R_d(s)$ for $(X_1,C_2')$ vanishes at $s=1$.  To 
obtain a
contradiction and prove the conjecture, we therefore need 
a nonvanishing
theorem which states that $R_d(1)$ is nonzero.  We shall 
in fact prove
that, if the genus of $C_2'$ is odd and its homology class 
is even
(conditions which eventually entail no loss of generality) 
and if the
conditions of Theorem~2 hold, then the value of $R_d(1)$ 
on the class
$\omega + \overline\omega$ of (c) is positive once $d$ is 
sufficiently large. 
If we recall again that the invariants for an immersed 
surface are defined
in terms of the embedded surface obtained by blowing up, 
we are led to blow
up $(X_1, C_2')$ at the single double point of $C_2'$ to 
obtain finally a
smooth algebraic pair $(X, C)$. The following result is 
therefore what is
wanted. 

\proclaim {Proposition 4}
Let $C$ be a smooth curve in an algebraic surface $X$\RM, 
satisfying the
conditions of Theorem~\RM2. Suppose in addition that $C$ 
has odd genus and that
its homology class has divisibility $2$ in $H_2(X, \Bbb 
Z)$.  Then for
large $d$ the value of the invariant $R_d(1)(\omega + 
\overline\omega)$ for the
pair $(X,C)$ is strictly positive. 
\endproclaim

Note that since $R_d(1)$ is just the sum of the invariants 
$q_{k,l}$ of a
given degree, it will suffice to show that each of these 
terms is
nonnegative and that at least one of them is positive.  
Under some
conditions on $\alpha$ and the metric on $X$, it was shown 
in \cite{5} that
the moduli spaces $M^\alpha_{k,l}$ can be interpreted as 
moduli spaces of
stable parabolic bundles. It is therefore tempting to try 
to adapt
Donaldson's argument in \cite{2} to prove that each 
$q_{k,l}$ is positive
when evaluated on the hyperplane class, provided that the 
degree $d$ is
large. (Actually this would not be quite the right thing; 
one should
construct a different version of $q_{k,l}$ which varies 
with $\alpha$, to
take account of how the natural polarization of a moduli 
space of parabolic
bundles changes as the parabolic weight is varied.)  
Unfortunately, there is
an obstruction to this programme.  A key technical step in 
the argument of
\cite{2} is to show that, once $d$ is large, the moduli 
spaces of stable
bundles on a complex surface have the dimension one would 
naively predict
from the index formula, namely, $d$.  This ``regularity'' 
result is false in
the context of parabolic bundles.

To explain what is true in the way of regularity, it is 
convenient to
introduce the ``magnetic charges'' $m_1 = k$ and $m_2 = k +
 l -
(\varSigma\cdot\varSigma)/4$ in place of $k$ and $l$.  The 
naively
expected complex dimension of the moduli space 
$M^\alpha_{k,l}$ is then $d
\sim 2m_1 + 2m_2$, according to the formula (1).  In order 
for the moduli
space to have this expected dimension, it is not enough 
that $d$ alone be
large---it is necessary for both $m_1$ and $m_2$ to be 
large. 

Fortunately, a general vanishing theorem was proved in 
\cite{5} which shows
that, if the absolute value of the difference $|m_1 - 
m_2|$ is larger
than a quantity $(K_X\cdot C)/4$, then the invariants
$q_{k,l}$ of a pair $(X,C)$ vanish when restricted to 
homology
classes orthogonal to $[C]$ in $H_2(X)$.  So, if we take 
such a
homology class, then the moduli spaces which might have 
the wrong dimension
(where only one of $m_1$ or $m_2$ is large) will not 
contribute.  Taking
this line forces us to abandon the idea of following 
Donaldson's argument
from \cite{2}, since the hyperplane class is not 
orthogonal to $C$.
Instead, we adapt O'Grady's argument from \cite{8}.

As was mentioned in the introduction, the results of 
\cite{8} show that,
under suitable conditions on a complex surface $X$ and 
$k$, the value
$q_k(\omega + \overline\omega)$ of the ordinary polynomial 
invariant is strictly
positive when $\omega$ is dual to a generic holomorphic 
$2$-form.  (Note
that such a class is always going to be orthogonal to a 
holomorphic curve
such as $C$.) Part of the argument adapts readily to the 
parabolic case to
show that $q_{k,l}(\omega + \overline\omega)$ is at least 
nonnegative, provided
only that the moduli spaces have the correct regularity 
properties.  All
that remains finally is to show that at least one of these 
values is
nonzero.  The last ingredient is another hard result from 
\cite{5} which
shows that, for the special value of the monopole number 
$l = (g-1)/2$, the
invariant $q_{k,l}$ for the pair $(X,C)$ is equal to $2^g 
q_k$ when
restricted to the orthogonal complement of $C$.  So, for 
this special value
of $l$, the nonvanishing of $q_{k,l}$ can be deduced from 
the nonvanishing of $q_k$.

\heading Acknowledgment\endheading

The author thanks Simon Donaldson, Bob Gompf, John Morgan, 
and
Kieran O'Grady for their help in preparing this paper, and,
in particular, Tom Mrowka for many hours of discussion, 
out of which
these results gradually emerged.


\enddocument